\documentclass[12pt]{amsart}
\usepackage{amsmath,amssymb, amsthm}

\theoremstyle{definition}

\newtheorem{lemma}{Lemma}
\newtheorem{theorem}{Theorem}

\newtheorem{corollary}{Corollary}
\def\fl#1{\left\lfloor#1\right\rfloor}

\title{Error terms for continued fractions of $e^{1/s}$ and $\sqrt{\frac{v}{u}}\tanh\!\Bigl(\frac{1}{\sqrt{uv}}\Bigr)$} 

\author{Nikita Kalinin}
\address{Mathematics and Computer Science Department,
Guangdong Technion Israel Institute of Technology,
Shantou 515603, China}
\email{nikita.kalinin@gtiit.edu.cn}

\author{Takao Komatsu}
\address{Institute of Mathematics,
Henan Academy of Sciences,
Zhengzhou 450046, China}
\address{Department of Mathematics,
Institute of Science Tokyo,
2-12-1 Ookayama, Meguro-ku,
Tokyo 152-8551, Japan}
\email{komatsu@zstu.edu.cn}
\email{komatsu.t.al@m.titech.ac.jp}

\date{}

\begin{document}
\maketitle

\begin{abstract}Many classical identities arise from nothing more mysterious than looking at the same object in two different ways. A number, a function, or a combinatorial object may admit several natural decompositions, and by disassembling it in one way and reassembling it in another, we often obtain unexpected corollaries. Telescoping sums provide a particularly vivid incarnation of this principle: by arranging terms so that successive contributions cancel, one performs a conceptual “cut-and-paste’’ that often admits a clean geometric interpretation. Generating functions offer a complementary perspective. Encoding a problem into a formal power series and then evaluating that series at a prescribed point naturally expresses the same quantity as an infinite (or finite) expansion, and equating these representations yields a wealth of identities.

For example, for a real number \(\alpha\) given by its continued fraction expansion 
$\alpha = [a_0, a_1,a_2,\dots]$, 
with convergents \(p_n/q_n\) and error terms 
$E_n := p_n - \alpha q_n$, 
one can obtain “additive'' decompositions of the form 
$\sum_{n\ge-1} a_{n+1}\,\lvert E_n\rvert \;=\; \alpha + 1$, 
$\sum_{n\ge-1} a_{n+1}\,E_n^{2} \;=\; \alpha$.  
Thus \(\alpha\) and \(\alpha+1\) themselves appear as weighted sums of the local approximation errors of their convergents.  In this note we explore what such decompositions yield in two explicit cases: the continued fraction
\[
e^{1/s} = [1;\,{\overline{(2k-1)s-1,1,1}}]_{k=1}^{\infty}
\]
and the continued fraction
\[
\frac{s}{u}\tanh\!\Bigl(\frac{1}{s}\Bigr)
  = [\,0;\,\overline{(4k-3)u,\,(4k-1)\tfrac{s^{2}}{u}}\,]_{k=1}^{\infty}
\]

\end{abstract}

\section{Introduction and background on error sums for continued fractions}

The continued fraction expansion of a real number $\alpha$ is one of the most classical and powerful tools in Diophantine approximation.  For an irrational $\alpha$, its simple continued fraction
\[
\alpha = [a_0; a_1,a_2,a_3,\dots]
\]
gives rise to convergents $p_n/q_n$ with excellent approximation properties:
\[
\left|\alpha-\frac{p_n}{q_n}\right| = \frac{|E_n|}{|q_n|}
\qquad\text{where}\qquad
E_n := p_n - \alpha q_n,
\]
and the \emph{error terms} $E_n$ encode subtle arithmetic information about $\alpha$.  Already in the eighteenth century, Euler and Lambert exploited continued fractions of special functions such as $e^x$ and $\tan x$ to prove irrationality results and to understand the quality of rational approximations to these values~\cite{euler1744fractionibus,lambert1761memoire,hermite1874fonction}.  In particular, the continued fraction for $e$ and its relatives has become a standard textbook example; see, for instance, classical accounts and expositions in~\cite{olds1970simple,matthews1970some,van1996continued}.

For certain special functions, the continued fraction expansion exhibits a striking regularity, with partial quotients given by arithmetic or geometric progressions.  Continued fractions of this \emph{Hurwitzian} type were studied systematically by Hurwitz and later Lehmer~\cite{hurwitz1887entwicklung,lehmer1918arithmetical}, who obtained structural results on their convergents and quantitative control on the approximation errors. By a Hurwitzian continued fraction we mean one whose partial quotients are eventually given by polynomial functions of $n$. In the twentieth and early twenty-first centuries, this line of work has been continued and refined in many directions: Tasoevian and other structured continued fractions, their Diophantine and combinatorial properties, and connections with hypergeometric functions have been investigated in, for example,~\cite{laughlin2019some,hetyei2014hurwitzian,komatsu2012some, tasoev2000rational}.

A particularly prominent family of examples comes from exponential and hyperbolic functions.  The continued fractions for $e^{1/M}$ and related exponential values have been derived systematically by van~der~Poorten and Osler, who gave modern derivations of the patterns in the partial quotients and discussed their arithmetic significance~\cite{van1996continued,osler2006proof}.  In a series of papers, Komatsu developed a systematic theory of convergent errors for such continued fractions: for $\alpha$ of the form $e^{1/s}$, $\tanh(1/s)$, or closely related special functions, he obtained exact integral representations for the errors $E_n$ 
often with a clean periodic dependence on $n$ modulo a fixed small integer, and in many cases these integrals admit further closed-form evaluation or hypergeometric simplification~\cite{komatsu2009diophantine,komatsu2009diophantine2,komatsu2012some}.  These results make the error terms completely explicit.

A complementary viewpoint, due to Ridley--Petruska and later developed extensively by Elsner and coauthors, is to study not individual $E_n$ but \emph{error sums} such as
\[
\sum_{n\ge0} (-E_n),
\qquad
\sum_{n\ge0} |E_n|,
\qquad
\sum_{n\ge0} E_n x^n.
\]
Ridley and Petruska introduced an ``error-sum function'' and expressed it in terms of the partial quotients of the continued fraction~\cite{ridley2000error}.  Elsner showed that for several important values, including $\alpha=e$, such sums can be evaluated in closed form and turn out to involve classical special functions such as the error function $\operatorname{erf}(x)$~\cite{elsner2011error,elsner2014error}.  These identities reveal global patterns in the sequence $(E_n)$ and connect Diophantine approximation with analytic properties of special functions.

\medskip

The goal of this note is to apply the following identities, established in
\cite{easter,certain},
\begin{align}
\sum_{n\ge -1} a_{n+1}\,|E_n| = \alpha + 1,
\label{eq11}
\\
\sum_{n\ge -1} a_{n+1}\,E_n^2 = \alpha,
\label{eq12}
\end{align}
where the coefficients $a_{n+1}$ are given explicitly by the continued--fraction expansion
$\alpha=[a_0;a_1,a_2,\dots]$.
These formulas express $\alpha$ itself as a weighted sum of its own approximation
errors $E_n=p_n-\alpha q_n$.
The identities are remarkable in that they recover $\alpha$ purely from the sequence
$\{E_n\}$, without reference to the convergents $(p_n,q_n)$ or to the continued fraction
algorithm generating them.
From this point of view, $(a_{n+1})$ plays the role of Fourier--type weights which extract
$\alpha$ from the full error profile.

Our main examples are 
\begin{itemize}
  \item the Hurwitzian continued fraction for $e^{1/s}=[1;\,\overline{(2k-1)s-1,1,1}]_{k=1}^{\infty}$,  
  \item and, for $s=\sqrt{u v}$, the continued fraction \cite{komatsu2009diophantine2}
  \[
  \alpha = \sqrt{\frac{v}{u}}\tanh\!\Bigl(\frac{1}{\sqrt{uv}}\Bigr)
  = \frac{s}{u}\tanh\!\Bigl(\frac{1}{s}\Bigr)
  = [\,0;\,\overline{(4k-3)u,\,(4k-1)v}\,]_{k\ge1}.
  \]
\end{itemize}

The paper is organized as follows.  In Section~\ref{sec:setting} we recall the continued fraction expansions for $e^{1/s}$ and $\tanh(1/s)$ and state the integral formulas for the corresponding error terms $E_n$.  Section~\ref{sec:series-expansions} develops the resulting factorial series in powers of $s^{-1}$ and applies the identities (\ref{eq11}) and (\ref{eq12}) to obtain two different expansions for the same quantities.  By comparing coefficients of $s^{-k}$, we derive several families of purely combinatorial identities, such as
\[
\sum_{n=0}^k (2n+1)\,\frac{k!}{(k-n)!(k+n+1)!} \;=\frac{1}{(k+1)!}\,
   {}_3F_2\!\left(\begin{matrix}
     -k,\ \tfrac32,\ 1\\[2pt]
     \tfrac12,\ k+2
   \end{matrix}\Biggm|\,-1\right)= \; \frac{1}{k!}
\]
\[
\frac{1}{k!}
\;=\;
\sum_{n=0}^{\lfloor k/2\rfloor}
(4n+1)\,
\frac{2}{(k-2n)!}\,
\frac{k!}{(2n+k+1)!},
\]
\[
\frac{1}{k!}
\;=\;
\sum_{n=0}^{\lfloor (k-1)/2\rfloor}
(4n+3)\,
\frac{2}{(k-2n-1)!}\,
\frac{k!}{(k+2n+2)!}.
\]

The last two identities also follow from the first one and the fact that  for every integer $k\ge 1$ we have
\[
{}_3F_2\!\left(\begin{matrix}
  -k,\ \tfrac32,\ 1\\[2pt]
  \tfrac12,\ k+2
\end{matrix}\Biggm|\,1\right)
:=\sum_{j=0}^k
\frac{(-k)_j\bigl(\tfrac32\bigr)_j(1)_j}{\bigl(\tfrac12\bigr)_j(k+2)_j\,j!}
=0.
\]

We also obtain a quadratic sum identity arising from the quadratic error-sum formula (\ref{eq12}) for $e^{1/s}$ and $\frac{s}{u}\tanh(1/s)$.

\section{Explicit error formulas for two Hurwitzian continued fractions}\label{sec:setting}

We now turn to the explicit computation of the error terms 
\(E_n=p_n-\alpha q_n\) for the two continued fractions described in the
introduction.  In both cases the structure of the partial quotients
leads to clean integral representations for \(E_n\), which in turn 
yield rapidly convergent factorial series.  In the sequel we will substitute these formulas into the
identities (\ref{eq11}) and (\ref{eq12}), obtaining two different series
expansions for the same quantity (namely \(\alpha\) or \(\alpha+1\)).
Since these series must agree term by term as formal power series in
\(s^{-1}\), equating the coefficients of \(s^{-k}\) produces a family of
new identities involving the factorial expressions appearing in the
error terms.  This section prepares the explicit expansions needed for
that coefficientwise comparison.

\subsection{Error Terms for $e^{1/s}$}

The continued fraction expansion of $\alpha=e^{1/s}$ is given by 
\[
\alpha=[1;\,\overline{(2k-1)s-1,1,1}]_{k=1}^{\infty}.
\]

Following the classical integral representation \cite{komatsu2009diophantine}, for every integer $n\ge0$,
\begin{align*}
E_{3n}
&=
-\frac{1}{s^{\,n+1}}
\int_{0}^{1}\frac{x^{n}(x-1)^{n}}{n!}\,e^{x/s}\,dx,
\\[2mm]
E_{3n+1}
&=
\frac{1}{s^{\,n+1}}
\int_{0}^{1}\frac{x^{n+1}(x-1)^{n}}{n!}\,e^{x/s}\,dx,
\\[2mm]
E_{3n+2}
&=
\frac{1}{s^{\,n+1}}
\int_{0}^{1}\frac{x^{n}(x-1)^{n+1}}{n!}\,e^{x/s}\,dx.
\end{align*}

Inserting the Taylor series 
\begin{equation}
e^{x/s}=\sum_{m\ge0}\frac{x^{m}}{m!\,s^{m}}
\label{eq:ser-e}
\end{equation}
and integrating termwise
\begin{equation}
\int_{0}^{1} x^{a}(1-x)^{b}\,dx=\frac{a!\,b!}{(a+b+1)!},
\label{eq:int}
\end{equation}
we obtain for every $n\ge0$ and $m\ge0$:
\begin{align*}
|E_{3n}| &= \sum_{m\ge0}\frac{(n+m)!}{m!\,(2n+m+1)!}\,s^{-(n+1+m)},\\[2mm]
|E_{3n+1}| &= \sum_{m\ge0}\frac{(n+1+m)!}{m!\,(2n+m+2)!}\,s^{-(n+1+m)},\\[2mm]
|E_{3n+2}| &= \sum_{m\ge0}\frac{(n+1)(n+m)!}{m!\,(2n+m+2)!}\,s^{-(n+1+m)}.
\end{align*}

\subsection{Error terms for $\sqrt{\frac{v}{u}}\;\tanh\!\Bigl(\frac{1}{s}\Bigr)$}

Let $u>0$ and $s>0$ with $s=\sqrt{u v}$, set 
\[
\alpha
=\sqrt{\frac{v}{u}}\;\tanh\!\Bigl(\frac{1}{s}\Bigr)
=\frac{s}{u}\,\tanh\!\Bigl(\frac{1}{s}\Bigr).
\]
Then $\alpha$ has the continued fraction expansion \cite{komatsu2009diophantine2}
\[
\alpha
=[\,0;\,\overline{(4k{-}3)u,\,(4k{-}1)v}\,]_{k=1}^{\infty}.
\]
Denote by $p_n/q_n$ the $n$th convergent, and set $D_n:=p_n-\alpha\,q_n$ (to distinguish it from $E_n$). Proceeding as in the exponential case, we insert the power series
expansion of $e^{2x/s}$ into the integrals and evaluate them termwise.
This again produces explicit factorial expressions for each $D_n$.

For $n\ge1$ the error terms satisfy the integral formulas (\cite{komatsu2009diophantine2}), which are obtained from the general formulas with $s^2=uv$ eliminated: 
\begin{align}
D_{2n-1}
&=
-\frac{1}{e^{2/s}+1}
\left(\frac{4}{s^2}\right)^{\!n}
\int_0^1 \frac{x^{2n-1}(x-1)^{2n-1}}{(2n-1)!}\,e^{2x/s}\,dx,
\label{d2n1}
\\[1mm]
D_{2n}
&=
-\frac{2}{(e^{2/s}+1)\,u}
\left(\frac{4}{s^2}\right)^{\!n}
\int_0^1 \frac{x^{2n}(x-1)^{2n}}{(2n)!}\,e^{2x/s}\,dx.
\label{d2n}
\end{align}
For $n\ge1$ these integrals converge absolutely.

Using 
\[
\frac{1}{e^{2/s}+1}
=\frac{1}{2}\bigl(1-\tanh(1/s)\bigr),
\]
we rewrite these as
\begin{align*}
D_{2n-1}
&=
\frac{1}{2}\bigl(1-\tanh(1/s)\bigr)\,
\left(\frac{4}{s^2}\right)^{\!n}
\int_0^1 \frac{x^{2n-1}(1-x)^{2n-1}}{(2n-1)!}\,e^{2x/s}\,dx,
\\[1mm]
D_{2n}
&=
-\frac{1-\tanh(1/s)}{u}\,
\left(\frac{4}{s^2}\right)^{\!n}
\int_0^1 \frac{x^{2n}(1-x)^{2n}}{(2n)!}\,e^{2x/s}\,dx.
\end{align*}

Similar to the exponent case, substituting into the integrals and cancelling the factorials in the denominators yields the convergent series (for each fixed $n\ge1$):
\begin{align*}
|D_{2n-1}|
&=
\frac{1}{2}\bigl(1-\tanh(1/s)\bigr)
\sum_{m\ge0}
\frac{4^n\,2^m}{s^{2n+m}\,m!}\,
\frac{(2n-1+m)!}{(4n+m-1)!},
\\[1mm]
|D_{2n}|
&=
(1-\tanh(1/s))\frac{1}{u}\,
\sum_{m\ge0}
\frac{4^n\,2^m}{s^{2n+m}\,m!}\,
\frac{(2n+m)!}{(4n+m+1)!}.
\end{align*}

\section{Weighted error-sum identities and combinatorial corollaries}\label{sec:series-expansions}
\subsection{The sum of weighted error terms for $e^{1/s}$}

In the general setting of \cite{easter,certain}, for a real number 
\(\alpha=[a_0;a_1,a_2,\dots]\) with convergents \(p_n/q_n\) and errors 
\(E_n=p_n-\alpha q_n\), there exist explicit weights \(a_{n+1}\) (depending on 
the partial quotients) such that (\ref{eq11}) and (\ref{eq12}) hold.  
We adopt the standard convention \(p_{-1}=1,q_{-1}=0\), so \(E_{-1}=1\),
and for our example \(\alpha=e^{1/s}\) we have \(a_0=1\), hence 
\(a_0|E_{-1}|=1\).

For \[
e^{1/s} = [1;\,{\overline{(2k-1)s-1,1,1}}]_{k=1}^{\infty}
\]
we have $a_{3n+1}=(2n+1)s-1$ and $a_{3n+2}=a_{3n+3}=1$ and the “\(n=-1\)” term contributes \(a_0|E_{-1}|=1\).
Let us expand both sides of (\ref{eq12}) in powers of $s^{-1}$ and compare coefficients.

Coefficients of $s^{-r}$ in $\sum_{n=-1}^\infty a_{n+1}E_n^2$ for $e^{1/s}$ are computed as follows:

\begin{theorem} For every integer $r\ge 0$, the coefficient of $s^{-r}$ in
$\sum_{n=-1}^{\infty} a_{n+1}E_n^2$ is given by
\[
\begin{aligned}
&\frac{1}{r!}=\bigl[s^{-r}\bigr]\!\Bigl(\sum_{n=-1}^{\infty} a_{n+1}\,E_n^2\Bigr)= \mathbf{1}_{\{r=0\}}
\\[1mm]
&\quad
 +\sum_{n=0}^{\infty}
\Biggl[
\;(2n+1)\sum_{\substack{i,j\ge0\\ 2n+1+i+j=r}}
\frac{(n+i)!}{i!\,(2n+i+1)!}\,\frac{(n+j)!}{j!\,(2n+j+1)!}
\\
&\quad
-\sum_{\substack{i,j\ge0\\ 2n+2+i+j=r}}
\frac{(n+i)!}{i!\,(2n+i+1)!}\,\frac{(n+j)!}{j!\,(2n+j+1)!}
\\
&\quad
+\sum_{\substack{i,j\ge0\\ 2n+2+i+j=r}}
\biggl(
\frac{(n+1+i)!}{i!\,(2n+i+2)!}\,\frac{(n+1+j)!}{j!\,(2n+j+2)!}\\
&\quad +
\frac{(n+1)\,(n+i)!}{i!\,(2n+i+2)!}\,
      \frac{(n+1)\,(n+j)!}{j!\,(2n+j+2)!}
\biggr)
\Biggr].
\end{aligned}
\]
\end{theorem}

We do not know a direct combinatorial or conceptual proof of this equality; at present
it is obtained indirectly by comparing the coefficient of \(s^{-r}\) in both
sides of (\ref{eq12}).  For each fixed \(r\), the sums over \(n,i,j\) are finite
(because of the constraints \(2n+1+i+j=r\) or \(2n+2+i+j=r\)), so the above
identity is purely combinatorial.

 Things are better for (\ref{eq11}), namely

\[
\sum_{n=-1}^\infty a_{n+1}|E_n|= \alpha+1 \;=\; 2+\sum_{k\ge1}\frac{1}{k!}\,s^{-k}.
\]
Thus the constant term is \(2\), and for each \(k\ge1\) the coefficient of \(s^{-k}\) in $\sum_{n=-1}^{\infty}a_{n+1}\,|E_n|$ is \(1/k!\).

To match this,  compute
\[
\bigl[s^{-k}\bigr]\!\Bigl(\sum_{n=-1}^{\infty} a_{n+1}\,|E_n|\Bigr)=\sum_{n=0}^{k}
(2n+1)\,\frac{k!}{(k-n)!\,\bigl(n+k+1\bigr)!}
\;+\;\]
\[
\sum_{n=0}^{k-1}
\biggl[
-\frac{(k-1)!}{(k-1-n)!\,\bigl(n+k\bigr)!}
\;+\;
\frac{k!}{(k-1-n)!\,\bigl(n+k+1\bigr)!}
\;+\;
\frac{(n+1)(k-1)!}{(k-1-n)!\,\bigl(n+k+1\bigr)!}
\biggr].
\]

Since 
$-(n+k+1)+k+(n+1)=0$ the second line of the above formula is identically zero.

Thus, the immediate corollary of our approach is as follows

\begin{corollary}
For every integer $k\ge 0$, one has
\begin{equation}
S_k:=\sum_{n=0}^k (2n+1)\frac{k!}{(k-n)!(n+k+1)!} = \frac{1}{k!}.
\label{eq1}
\end{equation}
\label{cor1}
\end{corollary}

For completeness, let us below derive this identity directly.

\subsection{Derivation of $\sum_{n=0}^k (2n+1)\frac{k!}{(k-n)!(n+k+1)!} = \frac{1}{k!}$}
We first rewrite the factorials in terms of binomial coefficients.
Recall that

\[
\frac{k!}{(k-n)!(k+n+1)!}
  =\frac{k!}{(2k+1)!}\binom{2k+1}{k-n}.
\]
Substituting this into the definition of \(S_k\) gives
\begin{equation*}
S_k
  =\frac{k!}{(2k+1)!}\sum_{n=0}^k (2n+1)\binom{2k+1}{k-n}.
\end{equation*}

Introducing a new index
$m:=k-n$, we rewrite 

\[
\sum_{n=0}^k (2n+1)\binom{2k+1}{k-n}
  =\sum_{m=0}^k (2k-2m+1)\binom{2k+1}{m}=\]
  \[=(2k+1)\sum_{m=0}^k\binom{2k+1}{m}
     -2\sum_{m=0}^k m\binom{2k+1}{m}.
\]

Next,
\[
\sum_{m=0}^k m\binom{2k+1}{m}
  =(2k+1)\sum_{m=0}^k\binom{2k}{m-1}
  =(2k+1)\sum_{r=0}^{k-1}\binom{2k}{r},
\]
where in the last step we put \(r=m-1\).

Finally, using Pascal's rule
\[
\binom{2k+1}{m}=\binom{2k}{m}+\binom{2k}{m-1}.
\]

and
\[
  \left(\sum_{m=0}^k\binom{2k}{m}
          +\sum_{r=0}^{k-1}\binom{2k}{r}\right)
     -2\sum_{r=0}^{k-1}\binom{2k}{r}
  =\sum_{m=0}^k\binom{2k}{m}
   -\sum_{r=0}^{k-1}\binom{2k}{r}
  =\binom{2k}{k},
\]

we obtain
\[
S_k=\frac{k!}{(2k+1)!}\,(2k+1)\binom{2k}{k}
   =\frac{k!}{(2k+1)!}\,\frac{(2k+1)!}{k!\,k!}
   =\frac{1}{k!},
\]
which is the desired identity.
\qed

\subsection{Hypergeometric reformulation of the identity for $S_k$}

The identity (\ref{eq1}) in Corollary \ref{cor1} 
also admits a compact hypergeometric reformulation.  We briefly record it here,
as it may be of independent interest.

Recall the Pochhammer symbol (rising factorial)
\[
(a)_n := a(a+1)\cdots(a+n-1)=\frac{\Gamma(a+n)}{\Gamma(a)}.
\]

First, express the linear factor \(2n+1\) via Pochhammer symbols with  half-integer parameters. Using the gamma function representation,
\[
\frac{(3/2)_n}{(1/2)_n}
 =\frac{\Gamma(\frac32+n)}{\Gamma(\frac32)}
  \frac{\Gamma(\frac12)}{\Gamma(\frac12+n)}
 = 2\,\frac{\Gamma(n+\frac32)}{\Gamma(n+\frac12)}
 = 2\Bigl(n+\tfrac12\Bigr)
 = 2n+1.
\]
Hence
\[
2n+1 = \frac{(3/2)_n}{(1/2)_n}.
\]

Next, we rewrite the factorial part using integer Pochhammer symbols.
For an integer \(k\ge0\) we have
\[
(-k)_n = (-1)^n \frac{k!}{(k-n)!},\qquad
(k+2)_n = \frac{(k+n+1)!}{(k+1)!}.
\]
Therefore
\[
\frac{k!}{(k-n)!(k+n+1)!}
 = \frac{1}{(k+1)!}\,(-1)^n\,
   \frac{(-k)_n}{(k+2)_n}.
\]

Substituting both identities into the summand of \(S_k\), we obtain
\[
(2n+1)\frac{k!}{(k-n)!(k+n+1)!}
  = \frac{(3/2)_n}{(1/2)_n}\,
    \frac{1}{(k+1)!}\,(-1)^n
    \frac{(-k)_n}{(k+2)_n}.
\]
Thus
\[
S_k
 = \frac{1}{(k+1)!}
   \sum_{n=0}^k
    (-1)^n
    \frac{(-k)_n\,(3/2)_n}{(1/2)_n\,(k+2)_n}.
\]

Now notice that \((1)_n = n!\). Insert the factor
\(\dfrac{(1)_n}{n!}=1\) into each term:
\begin{align*}
S_k
 &= \frac{1}{(k+1)!}
   \sum_{n=0}^k
   \frac{(-k)_n\,(3/2)_n\,(1)_n}{(1/2)_n\,(k+2)_n}\,
   \frac{(-1)^n}{n!}\\
&=\frac{1}{(k+1)!}\,
   {}_3F_2\!\left(\begin{matrix}
     -k,\ \tfrac32,\ 1\\[2pt]
     \tfrac12,\ k+2
   \end{matrix}\Biggm|\,-1\right),
\end{align*}
where 
\[
{}_3F_2(a_1,a_2,a_3; b_1,b_2; z)
 := \sum_{n=0}^{\infty}
    \frac{(a_1)_n(a_2)_n(a_3)_n}{(b_1)_n(b_2)_n}
    \frac{z^n}{n!} 
\]
is the generalized hypergeometric series.  

Since the parameter \(-k\) is a non-positive integer, the series
terminates at \(n=k\), so the infinite \({}_3F_2\) is actually a
finite sum and coincides with our original sum.  Using the identity
proved independently
\[
S_k=\frac{1}{k!},
\]
we equivalently have the closed hypergeometric form
\[
{}_3F_2\!\left(\begin{matrix}
     -k,\ \tfrac32,\ 1\\[2pt]
     \tfrac12,\ k+2
   \end{matrix}\Biggm|\,-1\right)
 = k+1.
\]

So, for each $k\geq 0$ we have the following identity.  

\begin{theorem}\label{thm3}
For every integer $k\ge 0$, one has
\[
\sum_{n=0}^k (2n+1)\frac{k!}{(k-n)!(k+n+1)!}
 = \frac{1}{(k+1)!}\,
   {}_3F_2\!\left(\begin{matrix}
     -k,\ \tfrac32,\ 1\\[2pt]
     \tfrac12,\ k+2
   \end{matrix}\Biggm|\,-1\right)
 = \frac{1}{k!}.
\]
\end{theorem}

\subsection{The sum of weighted error terms for $\frac{s}{u}\tanh\!\Bigl(\frac{1}{s}\Bigr)$}
We have
\[
\frac{s}{u}\tanh\!\Bigl(\frac{1}{s}\Bigr)
  = [\,0;\,\overline{(4k-3)u,\,(4k-1)\tfrac{s^2}{u}}\,]_{k\ge1}.
 \] 
So
\[
a_{2k-1}=(4k-3)\,u,\qquad
a_{2k}=(4k-1)\,\frac{s^2}{u}\qquad(k\ge1),
\]
and $a_0=0$ since $\alpha$ has zero integer part. Substituting $D_n$ into (\ref{eq11})
\[
\sum_{n=-1}^{\infty}a_{n+1}\,|D_n|\;=\;\alpha+1
\;=\;\frac{s}{u}\,\tanh\!\Bigl(\frac{1}{s}\Bigr)+1.
\]
yields an expression of the form
\[
\frac{s}{u}\,\tanh\!\Bigl(\frac{1}{s}\Bigr)+1
=
\bigl(1-\tanh(1/s)\bigr)\,
\Bigl(C_0(s)+\frac{s}{u}\,C_1(s)\Bigr),
\]
where the coefficient functions $C_0(s)$ and $C_1(s)$ are given by the double series
\begin{align*}
C_0(s)
&=
\sum_{n\ge0}\sum_{m\ge0}
(4n+1)\,
\frac{4^n\,2^m}{s^{2n+m}\,m!}\,
\frac{(2n+m)!}{(4n+m+1)!},
\\[1mm]
C_1(s)
&=
\sum_{n\ge0}\sum_{m\ge0}
(4n+3)\,
\frac{4^{\,n+1}\,2^m}{2\,s^{2n+1+m}\,m!}\,
\frac{(2n+1+m)!}{(4n+3+m)!}.
\end{align*}

Since $u$ is an independent parameter, we can compare the constant term and the coefficient of $s/u$ on both sides. This gives the two formal identities
\begin{align*}
\frac{1}{1-\tanh(1/s)}
&=
\sum_{n\ge0}\sum_{m\ge0}
(4n+1)\,
\frac{4^n\,2^m}{s^{2n+m}\,m!}\,
\frac{(2n+m)!}{(4n+m+1)!},
\\[1mm]
\frac{\tanh(1/s)}{1-\tanh(1/s)}
&=
\sum_{n\ge0}\sum_{m\ge0}
(4n+3)\,
\frac{4^{\,n+1}\,2^m}{2\,s^{2n+1+m}\,m!}\,
\frac{(2n+1+m)!}{(4n+3+m)!}.
\end{align*}

Let us see what we can get from the first identity
\begin{equation}\label{eq:goal-z}
\frac{1}{1-\tanh(1/s)}
=\frac{e^{2/s}+1}{2}=
\sum_{n\ge0}\sum_{m\ge0}
(4n+1)\,
\frac{4^n\,2^m}{m!}\,
\frac{(2n+m)!}{s^{2n+m}(4n+m+1)!}\,.
\end{equation}

The coefficient of $s^{-k}$ on the left-hand side of \eqref{eq:goal-z} is $1$ for $k=0$ and
\[
\frac{2^{k-1}}{k!}\quad(k\ge1).
\]

For $k=2n+m$, we compute the coefficient of $s^{-k}$ for the right-hand side as
\[
\sum_{n=0}^{\lfloor k/2\rfloor}
(4n+1)\,
\frac{2^k}{(k-2n)!}\,
\frac{k!}{(2n+k+1)!}.
\]

Thus, we get the following identity.

\begin{theorem}\label{thm4} For $k\ge1$
\[\frac{1}{k!}=\sum_{n=0}^{\lfloor k/2\rfloor}
(4n+1)\,
\frac{2}{(k-2n)!}\,
\frac{k!}{(2n+k+1)!},\]
or, in binomial form
\[(2k+1){\binom{2k}{k}}=2\sum_{n=0}^{\lfloor k/2\rfloor}(4n+1)\binom{2k+1}{k-2n}\]
\end{theorem}

We now perform a similar coefficient comparison for the second identity
\begin{multline*}
\sum_{k\ge1}\frac{2^{k-1}}{k!}\,s^{-k}=\frac{\tanh(1/s)}{1-\tanh(1/s)}=\frac{e^{2/s}-1}{2}\\
=
\sum_{n\ge0}\sum_{m\ge0}
(4n+3)\,
\frac{4^{\,n+1}\,2^m}{2\,s^{2n+1+m}\,m!}\,
\frac{(2n+1+m)!}{(4n+3+m)!}.
\end{multline*}

On the right-hand side, the general term is
\[
(4n+3)\,
\frac{4^{\,n+1}\,2^m}{2\,s^{2n+1+m}\,m!}\,
\frac{(2n+1+m)!}{(4n+3+m)!},
\]
so it contributes to the coefficient of $s^{-k}$ precisely when
\[
2n+1+m = k
\quad\Longleftrightarrow\quad
m = k-2n-1\ge 0.
\]
Hence, for fixed $k\ge1$, the admissible values of $n$ are
$0\le n\le \lfloor (k-1)/2\rfloor$, and the coefficient of $s^{-k}$ on
the right-hand side is
\[
\sum_{n=0}^{\lfloor (k-1)/2\rfloor}
(4n+3)\,
\frac{4^{\,n+1}\,2^{k-2n-1}}{2\,(k-2n-1)!}\,
\frac{(2n+1+k-2n-1)!}{(4n+3+k-2n-1)!}.
\]

Equating coefficients with $\sum_{k\ge1}\frac{2^{k-1}}{k!}\,s^{-k}$
gives the following theorem.

\begin{theorem}\label{thm5}
For every integer $k\ge1$,
\begin{equation}\label{eq:second-factorial-identity}
\frac{1}{k!}
=
\sum_{n=0}^{\lfloor (k-1)/2\rfloor}
(4n+3)\,
\frac{2}{(k-2n-1)!}\,
\frac{k!}{(k+2n+2)!}.
\end{equation}
Equivalently, multiplying both sides of \eqref{eq:second-factorial-identity}
by $(2k+2)!/k!$ and using
\[
\frac{(2k+2)!}{k!^2}
=(k+1)^2\binom{2k+2}{k+1},
\]
\[
\frac{(2k+2)!}{(k-2n-1)!(k+2n+2)!}
=(k+2n+3)\binom{2k+2}{k-2n-1},
\]
we obtain the binomial form
\[
(k+1)^2\binom{2k+2}{k+1}
=
2\sum_{n=0}^{\lfloor (k-1)/2\rfloor}
(4n+3)(k+2n+3)\,\binom{2k+2}{k-2n-1}.
\]
\end{theorem}

\noindent
As in the first case, this identity arises by comparing two expansions of
the same analytic function. One can deduce Theorem~\ref{thm4} and Theorem~\ref{thm5} from Theorem~\ref{thm3} by taking the sum and taking the difference with the hypergeometric series from the following lemma.

\begin{lemma}\label{lem:3F2-vanish}
For every integer $k\ge 1$ we have
\[
{}_3F_2\!\left(\begin{matrix}
  -k,\ \tfrac32,\ 1\\[2pt]
  \tfrac12,\ k+2
\end{matrix}\Biggm|\,1\right)
=\sum_{j=0}^k
\frac{(-k)_j\bigl(\tfrac32\bigr)_j(1)_j}{\bigl(\tfrac12\bigr)_j(k+2)_j\,j!}
=0.
\]
\end{lemma}

\subsection{Quadratic error-sum identity for $\frac{s}{u}\,\tanh\!\Bigl(\frac{1}{s}\Bigr)$}

From (\ref{d2n1}) and (\ref{d2n}), for $n\ge1$ we have
\begin{align*}
D_{2n-1}^2
&=\frac{1}{4}\bigl(1-\tanh(1/s)\bigr)^2\\
&\quad\times\sum_{m_1,m_2\ge0}
\frac{4^{2n}\,2^{m_1+m_2}}{s^{4n+m_1+m_2}\,m_1!\,m_2!}\,
\frac{(2n-1+m_1)!(2n-1+m_2)!}{(4n+m_1-1)!(4n+m_2-1)!},
\end{align*}
and
\begin{align*}
D_{2n}^2
&=\bigl(1-\tanh(1/s)\bigr)^2\frac{1}{u^2}\\
&\quad\times\sum_{m_1,m_2\ge0}
\frac{4^{2n}\,2^{m_1+m_2}}{s^{4n+m_1+m_2}\,m_1!\,m_2!}\,
\frac{(2n+m_1)!(2n+m_2)!}{(4n+m_1+1)!(4n+m_2+1)!}.
\end{align*}

Using $a_{2n+1}=(4n+1)u$ and $a_{2n}=(4n-1)\dfrac{s^2}{u}$, we compute the coefficient of $s^{-k}$ in (\ref{eq12})
\[
\sum_{n=-1}^{\infty}a_{n+1}\,D_n^2
=
\sum_{n\ge0}a_{2n+1}\,D_{2n}^2
+
\sum_{n\ge1}a_{2n}\,D_{2n-1}^2.
\]
Since $a_0=0$, there is no $n=-1$ contribution.  Substituting the above expressions for
$D_{2n}^2$ and $D_{2n-1}^2$ and collecting the common factor
$\bigl(1-\tanh(1/s)\bigr)^2/u$, we obtain
\[
\sum_{n=-1}^{\infty}a_{n+1}\,D_n^2
=
\frac{\bigl(1-\tanh(1/s)\bigr)^2}{u}
\Bigl(\Phi_0(s)+\Phi_1(s)\Bigr),
\]
where
\begin{align*}
\Phi_0(s)
&:=
\sum_{n\ge0}\sum_{m_1,m_2\ge0}
(4n+1)\\
&\quad\times
\frac{4^{2n}\,2^{m_1+m_2}}{s^{4n+m_1+m_2}\,m_1!\,m_2!}\,
\frac{(2n+m_1)!(2n+m_2)!}{(4n+m_1+1)!(4n+m_2+1)!},\\[1mm]
\Phi_1(s)
&:=
\sum_{n\ge1}\sum_{m_1,m_2\ge0}
(4n-1)\\
&\quad\times\frac{4^{2n-1}\,2^{m_1+m_2}}{s^{4n+m_1+m_2-2}\,m_1!\,m_2!}\,
\frac{(2n-1+m_1)!(2n-1+m_2)!}{(4n+m_1-1)!(4n+m_2-1)!}.
\end{align*}
For each fixed $k$, the constraints on $4n+m_1+m_2$ show that the
coefficient of $s^{-k}$ in $\Phi_0(s)+\Phi_1(s)$ is a finite sum.

By the quadratic identity (\ref{eq12}) we also have
\[
\sum_{n=-1}^{\infty}a_{n+1}\,D_n^2
=\alpha
=\frac{s}{u}\,\tanh\!\Bigl(\frac{1}{s}\Bigr).
\]
Cancelling the factor $1/u$ and dividing by $(1-\tanh(1/s))^2$ gives
\begin{equation}\label{eq:quad-master}
\Phi_0(s)+\Phi_1(s)
=\frac{s\,\tanh(1/s)}{\bigl(1-\tanh(1/s)\bigr)^2} = s\frac{e^{4/s}-1}{4}.
\end{equation}

Hence we obtain the
explicit exponential form
\begin{equation}\label{eq:quad-exp}
\Phi_0(s)+\Phi_1(s)
=\frac{s}{4}\Bigl(e^{4/s}-1\Bigr)
=\sum_{k\ge0}\frac{4^k}{(k+1)!}\,s^{-k}.
\end{equation}

On the other hand, the coefficient of $s^{-k}$ in $\Phi_0(s)$ comes from
all triples $(n,m_1,m_2)$ with $4n+m_1+m_2=k$, while the coefficient
in $\Phi_1(s)$ comes from triples with $4n+m_1+m_2=k+2$.  Thus, for
each integer $k\ge0$, equating coefficients in \eqref{eq:quad-exp}
gives the following identity.
%

\begin{theorem}
For every integer $k\ge0$,
\[
\frac{4^k}{(k+1)!}
=
\sum_{\substack{n\ge0,\ m_1,m_2\ge0\\ 4n+m_1+m_2=k}}
(4n+1)\,
\frac{4^{2n}\,2^{m_1+m_2}}{m_1!\,m_2!}\,
\frac{(2n+m_1)!(2n+m_2)!}{(4n+m_1+1)!(4n+m_2+1)!}
+\]
\[+
\sum_{\substack{n\ge1,\ m_1,m_2\ge0\\ 4n+m_1+m_2=k+2}}
(4n-1)\,
\frac{4^{2n-1}\,2^{m_1+m_2}}{m_1!\,m_2!}\,
\frac{(2n-1+m_1)!(2n-1+m_2)!}{(4n+m_1-1)!(4n+m_2-1)!}.
\]
For each fixed $k$ the two sums on the right-hand side are finite.
\end{theorem}

\section{Error terms for $e^{2/s}$}   

There are many ways to obtain new and also known identities by using (\ref{eq11}) or (\ref{eq12}).   
In \cite{komatsu2007integers} (see also \cite{komatsu2009diophantine2}), for odd integers $s\ge 3$, the error terms of $e^{2/s}$, whose continued fraction expansion is given by 
$$
e^{2/s}=[1;\overline{\frac{(6 k-5)s-1}{2},(12 k-6)s,\frac{(6 k-1)s-1}{2},1,1}]_{k=1}^\infty\,,
$$ 
were found explicitly:  
\begin{align*}
E_{5 n}&=-\left(\frac{2}{s}\right)^{3 n+1}\int_0^1\frac{x^{3 n}(x-1)^{3 n}}{(3 n)!}e^{2 x/s}d x\,,\\
E_{5 n+1}&=-\frac{2^{3 n+1}}{s^{3 n+2}}\int_0^1\frac{x^{3 n+1}(x-1)^{3 n+1}}{(3 n+1)!}e^{2 x/s}d x\,,\\
E_{5 n+2}&=-\left(\frac{2}{s}\right)^{3 n+3}\int_0^1\frac{x^{3 n+2}(x-1)^{3 n+2}}{(3 n+2)!}e^{2 x/s}d x\,,\\
E_{5 n+3}&=\left(\frac{2}{s}\right)^{3 n+3}\int_0^1\frac{x^{3 n+3}(x-1)^{3 n+2}}{(3 n+3)!}e^{2 x/s}d x\,,\\
E_{5 n+4}&=\left(\frac{2}{s}\right)^{3 n+3}\int_0^1\frac{x^{3 n+2}(x-1)^{3 n+3}}{(3 n+2)!}e^{2 x/s}d x\,. 
\end{align*} 
Then by applying (\ref{eq:ser-e}) and (\ref{eq:int}), we get  
\begin{align*}  
E_{5 n}&=p_{5 n}-q_{5 n}e^{2/s}\\
&=-\left(\frac{2}{s}\right)^{3 n+1}\int_0^1\frac{x^{3 n}(x-1)^{3 n}}{(3 n)!}\sum_{m=0}^\infty\frac{(2 x)^m}{m!s^m}d x\\
&=-\left(\frac{2}{s}\right)^{3 n+1}\sum_{m=0}^\infty\frac{2^m}{m!s^m}\int_0^1\frac{x^{m+3 n}(1-x)^{3 n}(-1)^n}{(3 n)!}d x\\ 
&=-\left(\frac{2}{s}\right)^{3 n+1}\sum_{m=0}^\infty\frac{2^m}{m!s^m}\frac{(-1)^n(m+3 n)!(3 n)!}{(3 n)!(m+6 n+1)!}\\
&=(-1)^{n+1}\sum_{m=0}^\infty\frac{(m+3 n)!}{m!(m+6 n+1)!}\left(\frac{2}{s}\right)^{m+3 n+1}\,.  
\end{align*} 
Hence, we have 
$$
|E_{5 n}|=\sum_{m=0}^\infty\frac{(m+3 n)!}{m!(m+6 n+1)!}\left(\frac{2}{s}\right)^{m+3 n+1}\,.   
$$ 

Similarly, we have  
\begin{align*}  
|E_{5 n+1}|&=\sum_{m=0}^\infty\frac{(m+3 n+1)!}{m!(m+6 n+3)!}\frac{2^{m+3 n+1}}{s^{m+3 n+2}}\,,\\ 
|E_{5 n+2}|&=\sum_{m=0}^\infty\frac{(m+3 n+2)!}{m!(m+6 n+5)!}\left(\frac{2}{s}\right)^{m+3 n+3}\,,\\ 
|E_{5 n+3}|&=\sum_{m=0}^\infty\frac{(m+3 n+3)!}{m!(m+6 n+6)!(3 n+3)}\left(\frac{2}{s}\right)^{m+3 n+3}\,,\\ 
|E_{5 n+4}|&=\sum_{m=0}^\infty\frac{(m+3 n)!(3 n+3)}{m!(m+6 n+1)!}\left(\frac{2}{s}\right)^{m+3 n+3}\,.   
\end{align*}

Since for $n\ge 0$ 
\begin{align*}
&a_{5 n+1}=\frac{(6 n+1)s-1}{2},\, a_{5 n+2}=(12 n+6)s,\,a_{5 n+3}=\frac{(6 n+5)s-1}{2},\\
&a_{5 n+4}=a_{5 n+5}=1\,, 
\end{align*} 
we have 

\begin{align*}  
&e^{2/s}+1=\sum_{n=-1}^\infty a_{n+1}|E_{n}|=1+\sum_{n=0}^\infty a_{n+1}|E_{n}|\quad(E_{-1}=1)\\
&=1+\sum_{n=0}^\infty\sum_{j=0}^4 a_{5 n+j+1}|E_{5 n+j}|\\
&=1+\sum_{n=0}^\infty \frac{(6 n+1)s-1}{2}\sum_{m=0}^\infty\frac{(m+3 n)!}{m!(m+6 n+1)!}\left(\frac{2}{s}\right)^{m+3 n+1}\\
&\quad +\sum_{n=0}^\infty (12 n+6)s\sum_{m=0}^\infty\frac{(m+3 n+1)!}{m!(m+6 n+3)!}\frac{2^{m+3 n+1}}{s^{m+3 n+2}}\\
&\quad +\sum_{n=0}^\infty \frac{(6 n+5)s-1}{2}\sum_{m=0}^\infty\frac{(m+3 n+2)!}{m!(m+6 n+5)!}\left(\frac{2}{s}\right)^{m+3 n+3}\\
&\quad +\sum_{n=0}^\infty \sum_{m=0}^\infty\frac{(m+3 n+3)!}{m!(m+6 n+6)!(3 n+3)}\left(\frac{2}{s}\right)^{m+3 n+3}\\
&\quad +\sum_{n=0}^\infty \sum_{m=0}^\infty\frac{(m+3 n+2)!(3 n+3)}{m!(m+6 n+6)!}\left(\frac{2}{s}\right)^{m+3 n+3}\\
\end{align*}
\begin{align*} 
&=1+\sum_{n,m\ge 0}\frac{(6 n+1)(m+3 n)!}{m!(m+6 n+1)!}\left(\frac{2}{s}\right)^{m+3 n}-\sum_{n,m\ge 0}\frac{(m+3 n)!}{m!(m+6 n+1)!}\frac{2^{m+3 n}}{s^{m+3 n+1}}\\
&\quad +\sum_{n,m\ge 0}\frac{(12 n+6)(m+3 n+1)!}{m!(m+6 n+3)!}\frac{2^{m+3 n+1}}{s^{m+3 n+1}}\\
&\quad 
+\sum_{n,m\ge 0}\frac{(6 n+5)(m+3 n+2)!}{m!(m+6 n+5)!}\left(\frac{2}{s}\right)^{m+3 n+2}\\
&\quad -\sum_{n,m\ge 0}\frac{(m+3 n+2)!}{m!(m+6 n+5)!}\frac{2^{m+3 n+2}}{s^{m+3 n+3}}\\
&\quad +\sum_{n,m\ge 0}\frac{(m+3 n+3)!}{3 m!(m+6 n+6)!(n+1)}\left(\frac{2}{s}\right)^{m+3 n+3}\\
&\quad +\sum_{n,m\ge 0}\frac{3(m+3 n+2)!(n+1)}{m!(m+6 n+6)!}\left(\frac{2}{s}\right)^{m+3 n+3}\,. 
\end{align*}

Comparing the coefficients with 
$$
e^{2/s}+1=2+\sum_{\ell=1}^\infty\frac{1}{\ell!}\left(\frac{2}{s}\right)^{\ell}\,,
$$ 
we have the following identity.  

\begin{theorem}  
For a positive integer $\ell$,  
\begin{align*}  
\frac{2^\ell}{\ell!}&=\sum_{n=0}^{\fl{\ell/3}}\frac{3\cdot 2^{\ell+1}(n+1)\ell!}{(\ell-3 n)!(\ell+3 n)!}-\sum_{n=0}^{\fl{(\ell-1)/3}}\frac{2^{\ell-1}(\ell-1)!}{(\ell-3 n-1)!(\ell+3 n)!}\\
&\quad +\sum_{n=0}^{\fl{(\ell-1)/3}}\frac{3\cdot 2^{\ell+1}(2 n+1)\ell!}{(\ell-3 n-1)!(\ell+3 n+2)!}\\ 
&\quad +\sum_{n=0}^{\fl{(\ell-2)/3}}\frac{2^{\ell}(6 n+5)\ell!}{(\ell-3 n-2)!(\ell+3 n+3)!} -\sum_{n=0}^{\fl{(\ell-3)/3}}\frac{2^{\ell-1}(\ell-1)!}{(\ell-3 n-3)!(\ell+3 n+2)!}\\ 
&\quad +\sum_{n=0}^{\fl{(\ell-3)/3}}\frac{2^{\ell}\ell!}{3(n+1)(\ell-3 n-3)!(\ell+3 n+3)!}+\sum_{n=0}^{\fl{(\ell-3)/3}}\frac{3\cdot 2^{\ell}(n+1)(\ell-1)!}{(\ell-3 n-3)!(\ell+3 n+3)!}\,. 
\end{align*}
\end{theorem}

\section*{Final comments}

We do not know a direct combinatorial or conceptual proof of the 
identities in this article, though though finding one would be interesting. 

\section*{Acknowledgement}  
The authors thank the referees for useful comments.  
The work of T.K. was partly supported by JSPS KAKENHI Grant Number 24K22835.


\end{document}